\documentclass[onecolumn]{IEEEtran}
\usepackage{amsmath,amsfonts}
\usepackage{algorithmic}
\usepackage{algorithm}
\usepackage{array}
\usepackage[caption=false,font=normalsize,labelfont=sf,textfont=sf]{subfig}
\usepackage{textcomp}
\usepackage{stfloats}
\usepackage{url}
\usepackage{verbatim}
\usepackage{graphicx}
\usepackage{cite}
\begin{document}

\title{New Quantum LDPC Codes Based on Euclidean Geometry}
\author{{Ya'nan Feng, Chuchen Tang, Chenming Bai. }
\thanks{This work is funded by Science and Technology Project of Hebei Education Department(Grant No. QN2020196) , Natural Science Foundation of Hebei Province(Grant No. A2021210027).}
\thanks{Ya'nan Feng(corresponding author), Chuchen Tang and Chenming Bai are with the Department of Mathematics and Physics, Shijiazhuang Tiedao University, Shijiazhuang, 050043, PR China(e-mail: fengfei-03@163.com;1202262012@student.stdu.edu.cn;baichm@stdu.edu.cn.).}
\thanks{Manuscript received July 25, 2023.}}



\maketitle

\begin{abstract}
  With the development of quantum error correction techniques, quantum low density parity check (QLDPC) codes become a promising area in quantum error correction codes. In this paper, the requirements of QLDPC codes based on points except the origin and lines not passing through the origin of Euclidean geometry are given. QLDPC codes based on all the lines and parallel classes are obtained respectively.
\end{abstract}


\section{Introduction}
\IEEEPARstart{L}{ow} density parity check (LDPC) codes are a class of classical codes that were first described by Gallager in 1962[1]. LDPC codes were rediscovered and presented in a graphical interpretation by Tanner in 1981[2]. There have been several notable attempts to construct regular and irregular good LDPC codes using algebraic combinatorics and random constructions. Liva et al.[3] presented a survey of the previous work done on algebraic constructions of LDPC codes based on finite geometries, elements of finite fields and RS codes.

Quantum information is sensitive to noise and needs error correction, control and recovery strategies. Quantum error correction codes are means to protect quantum information against noise and decoherence. But they have the drawback of requiring a growing number of physical qubits per logical qubit for error measurement. LDPC codes solved similar issues in classical coding theory. Therefore, it was natural to consider quantum low density parity check(QLDPC) codes. QLDPC codes with constant encoding rate can reduce overhead in fault-tolerant quantum computation to a constant level and had the potential to improve the efficiency of quantum computation. Besides, QLDPC codes were not just about fast decoding schemes, but might be the way towards a better understanding of quantum channels from a purely information theory viewpoint.

The first example of QLDPC code was conceived by Postol[4] in 2001, which was a Calderbank-Shor-Steane(CSS)-based non-dual-containing QLDPC code from a small finite geometry. Z. Babar et al.[5] gave a survey of QLDPC codes from 2000 to 2015, and classified them into four classes: dual-containing CSS codes, non-dual containing CSS codes, non-CSS codes and entanglement-assisted codes. There was a mount of work about QLDPC codes from finite geometry. Based on classical LDPC codes using the unicycle code design, Aly[6] gave the dual-containing QLDPC codes from finite geometry, especially from Euclidean geometry in 2008. Cao et al.[7] proposed a novel class of QLDPC codes constructed from cyclic classes of lines in Euclidean geometry in 2012. Farinholt[8] presented several new classes of QLDPC codes using finite projective planes in 2012. Fu[9] constructed asymmetric QLDPC codes using classical quasi-cyclic LDPC codes based on Euclidean geometry in 2018. Popatia[10] gave an explicit construction of a non CSS QLDPC code from projective geometry in 2022.

In this paper, we construct QLDPC codes based on Euclidean geometry. Section II introduces necessary notations and concepts. Section III gives the requirements of QLDPC codes constructed in [6] and presents a new construction method for the cases not included in the conditions. Section IV is devoted to QLDPC codes based on all the points and lines of Euclidean geometry.

\section{EUCLIDEAN GEOMETRY AND QUANTUM CODES}
\subsection{Euclidean Geometry $EG(m,q)$}
Let $q$ be a power of a prime $p$, i.e. $q=p^s$ and $s$ a positive integer. Let $EG(m,q)$ be the $ m $-dimensional Euclidean geometry over $\textbf{F}_q$ for some integer $m\geq 2$, see details in [11]. There are $q^m$ points in $EG(m,q)$ and every point is represented by an $ m $-tuple over $\textbf{F}_q$. There are $q^{m-1}(q^m-1)/(q-1)$ lines and each line can be described by a 1-dimensional subspace of the vector space $\textbf{F}_q^m$ or a coset of it. These lines can be partitioned into $(q^m-1))(q-1)$ parallel classes. For any point in $EG(m,q)$, there are $(q^m-1)/(q-1)$ lines intersect at this point. Any two points can define one and only one line in between. Each line passes through $q$ points and has $q^{m-1}-1$ lines parallel to it. Two lines intersect in at most one point. \\

\textbf{Example 1.}  Let $m=q=2$. There are 4 points and 6 lines in $EG(2,2)$. They are as follows.\\
$$
\begin{aligned}
	v_1=(0,0),\
	v_2=(0,1),\
	v_3=(1,0),\
	v_4=(1,1).\
\end{aligned}$$
$$
\begin{aligned}
	l_1=L(0,1)=\{(0,0),(0,1)\},\
	l_4=L(0,1)+(1,0)=\{(1,0),(1,1)\},\\
    l_2=L(1,0)=\{(0,0),(1,0)\},\
    l_5=L(1,0)+(0,1)=\{(0,1),(1,1)\},\\
    l_3=L(1,1)=\{(0,0),(1,1)\},\
	l_6=L(1,1)+(0,1)=\{(0,1),(1,0)\}.
\end{aligned}
$$

\subsection{Quantum Codes}
\textbf{Definition 1}\textbf{.} A $(\rho,\lambda)$ regular LDPC code is defined by a sparse binary parity check matrix $\textbf{H}$ satisfying the following properties.

i) $\rho$ is the number of ones in a column.

ii) $\lambda$ is the number of ones in a row.

iii) Any two rows have at most one nonzero element in common. The code does not have cycles of length four in its Tanner graph.

iv) $\rho$ and $\lambda$ are small in comparison to the number of rows and length of the code. In addition, rows of the matrix $\textbf{H}$ are not necessarily linearly independent.\\

Let $P=\{I, X, Z, Y\}$ be a set of Pauli matrices, where
$$
I=\begin{pmatrix}	1 & 0 \\0 & 1  \end{pmatrix},
X=\begin{pmatrix}	0 & 1 \\1 & 0  \end{pmatrix},
Y=\begin{pmatrix}	0 & -i \\i & 0  \end{pmatrix},
Z=\begin{pmatrix}	1 & 0 \\0 & -1  \end{pmatrix}.
$$
Most known quantum codes are stabilizer codes. Let $$S_j=E_1\otimes \cdots \otimes E_n,\; E_i \in P,\; j = 1, 2, \cdots, n-k. $$
Then $S_j$ can be seen as a binary vector of length $2n$.

Assume we have a stabilizer group $\mathcal{S}$ generated by the set
$$S=\{S_1, \cdots, S_{n-k}\},$$
such that every two operators commute with each other. A quantum code $Q$ is defined as the $+1$ joint eigenstates of stabilizer $\mathcal{S}$. That is, a codeword state $|\psi\rangle$ belongs to the code $Q$ if and only if
$$S_j|\psi\rangle=|\psi\rangle, \;\text{for all $S_j \in S$}. $$

The most common family of quantum stabilizer codes are CSS codes invented independently by Calderbank and Shor[12], as well as Steane[13,14]. If we take two binary classical codes $C_1$: $[n, k_1, d_1]$ and $C_2$: $[n, k_2, d_2]$ with parity check matrices $\textbf{H}_1$, $\textbf{H}_2$ respectively, satisfying that $C_1^{\perp} \subset C_2$, i.e. $\textbf{H}_1\textbf{H}_2^T=\textbf{O}$,  then the quantum code $Q$ with the stabilizer
$$
\textbf{S}=\bordermatrix{
	&n &n\cr
	&\textbf{H}_1 &\textbf{O} \cr
	&\textbf{O} &\textbf{H}_2 \cr
}\begin{array}{l}
	n-k_1\\n-k_2  \end{array}
$$
is a CSS code with parameters $[[n, k_1+k_2-n, d]]$, where each row of $\textbf{S}$ corresponds to a stabilizer of $\mathcal{S}$, $k_1+k_2>n$ and
$d=min\{d_1, d_2\}$. A particularly nice scenario is one in which a classical $[n, k, d]$ code $C$ with self-orthogonal parity check matrix $\textbf{H}(\textbf{HH}^T=\textbf{O})$, and the stabilizer becomes the following form
$$
\textbf{S}=\bordermatrix{
	&n &n\cr
	&\textbf{H} &\textbf{O} \cr
	&\textbf{O} &\textbf{H} \cr
}\begin{array}{l}
	n-k\\n-k  \end{array}.
$$
The quantum code constructed in this way has parameters $[[n, 2k-n, d]]$.

\section{QLDPC CODES BASED ON PARTIAL POINTS AND LINES OF $EG(m,q)$}
In this section, we mainly give the conditions of the QLDPC codes constructed in [6] and present a new construction method for the cases not included in the conditions,
\subsection{Requirements of the QLDPC codes in [6]}
Aly presented a construction of QLDPC codes based on the points and lines of Euclidean geometry, but the conditions were not given. We will show the requirements of the QLDPC codes proposed in [6].

Consider the points excluding the original point $\textbf{0}$ and lines not passing through $\textbf{0}$ in $EG(m,q)$. We can define the binary matrix $\textbf{H}_1=(h^{(1)}_{i,j})$, whose rows are indexed by lines and columns are indexed by points. $h^{(1)}_{i,j}=1$ if line $l_i$ passes through the point $v_j$, otherwise, $h^{(1)}_{i,j}=0$. Take $\textbf{H}_1^T$ to be the transpose of $\textbf{H}_1$, then $\textbf{H}_1^T$ has the following properties:

i) The total number of columns is  $(q^{m-1}-1)(q^m-1)/(q-1)$.

ii) The number of rows is given by $q^m-1$.

iii) $\lambda  = q(q^{m-1}-1)/(q-1)$ is the row weight.

iv) $\rho=q$ is the column weight.

v) Any two rows in $\textbf{H}_1^T$ have exactly one nonzero element in common. Similarly, any two columns have at most one nonzero element in common.

Next, we can define a self-orthogonal parity check matrix $\textbf{H}_1^{orth}$ based on $\textbf{H}_1^T$. Let $\textbf{1}$ be the $(q^m-1) \times 1$ column vector defined as $\textbf{1} = (1, \cdots, 1)^T$. If the row weight of $\textbf{H}_1^T$ is odd, then we can add the vector $\textbf{1}$ to form the matrix $\textbf{H}_1^{orth} = [\;\textbf{H}_1^T\;|\;\textbf{1}\;]$. Also, if the weight of a row in $\textbf{H}_1^T$ is even, then we can add the vector $\textbf{1}$ along with the $(q^m-1) \times (q^m-1)$ identity matrix to form the matrix $\textbf{H}_1^{orth} = [\;\textbf{H}_1^T\;|\;\textbf{1}\;|\;\textbf{I}\;]$.

Now the row weight of $\textbf{H}_1^T$ is $$\lambda=q(q^{m-1}-1)/(q-1)=q+q^2+\cdots+q^{m-1}.$$ $\lambda$ is odd only when $q$ is odd and $m$ is even. On the other hand, $\lambda$ is even if and only if that $q$ is even, or both $q$ and $m$ are odd.
Therefore, we have
$$
\textbf{H}_1^{orth} = \left\{
\begin{aligned}
\;&
\begin{split}
	[ \;\textbf{H}_1^T\;|\;\textbf{1}\;],               \ \ \;\text{if $q$ is odd, $m$ is even with $m\geq 4$}.\end{split} \\
	[& \;\textbf{H}_1^T\;|\;\textbf{1}\;|\;\textbf{I}\;],\;\text{if $q$ is even, or both $q$ and $m$ are odd}.
\end{aligned}
\right.
$$
For brevity, let $B_t=\frac{q^{t}-1}{q-1}$. It will be used frequently in the sequel.\\

\textbf{Lemma 2.} $\textbf{H}_1^{orth}$ is self-orthogonal.

\textit{Proof.} We need to show that $\textbf{H}_1^{orth} \cdot (\textbf{H}_1^{orth})^T=\textbf{O}$, that is, every row of $\textbf{H}_1^{orth}$ has even number of ones and any two distinct rows has even number of ones in common.

Any two different columns of $\textbf{H}_1$ has just one ``1'' in common, therefore, any two distinct rows of $\textbf{H}_1^{orth}$ intersect in exactly two nonzero positions.

(i) If $q$ is odd, $m$ is even with $m\geq4$, the row weight of $\textbf{H}_1^{orth}$ is
$$qB_{m-1}+1=1+q+\cdots+q^{m-1}.$$
Obviously, the row weight is an even number.

(ii) When $q$ is even, or both $q$ and $m$ are odd, the row weight of $\textbf{H}_1^{orth}$ is
$$qB_{m-1}+2=2+q+\cdots+q^{m-1},$$
which is clearly even. \\

\textbf{Theorem 3.} The QLDPC code $Q$ with stabilizer $\textbf{S}$ has parameters $[[n_1,k_1,\geq d_1]]$, where
$$
\textbf{S}=
	\begin{pmatrix}
		\textbf{H}_1^{orth} &\textbf{O}\\
		\textbf{O} &\textbf{H}_1^{orth}
	\end{pmatrix},
$$
and
$(n_1,k_1, d_1)=$
$$
\left\{
\begin{aligned}
&\begin{split}\;
((q^m-1)B_{m-1}+1,(q^m-1)(B_{m-1}-2)+1,q+1),\ \ \ \ \ \ \ \ \text{if $q$ is odd, $m$ is even with $m\geq 4$.}
\\
\end{split}\\~\\
&\begin{split}\;
((q^m-1)(B_{m-1}+1)+1,(q^m-1)(B_{m-1}-1)+1,q+1),\ \text{if $q$ is even, or both $q$ and $m$ are odd.}
	\end{split}
\end{aligned}
\right.
$$

\textit{Proof.} (i) If $q$ is odd and $m$ is even with $m\geq 4$,
$$
\begin{aligned}
	\textbf{S}=
	&\bordermatrix{
		&(q^m-1)B_{m-1}+1&(q^m-1)B_{m-1}+1\cr
		&\textbf{H}_1^{orth} & \textbf{O} \cr
		&\textbf{O} &\textbf{H}_1^{orth} \cr
	}\begin{array}{l}
		q^m-1\\q^m-1  \end{array}.
\end{aligned}
$$
It is apparent that the QLDPC code Q with stabilizer S has codelength $(q^m-1)B_{m-1}+1.$ The dimension of Q is
$$(q^m-1)B_{m-1}+1-2(q^m-1)=\frac{(q^m-1)(q^{m-1}-2q+1)}{q-1}+1,$$
which is positive when $m\geq4$.
Since the LDPC code having $\textbf{H}_1^{orth}$ as its parity check matrix has minimum distance at least $q+1,$ so it is with $Q$.

(ii) If $q$ is even, or both $q$ and $m$ are odd, the number of columns of $\textbf{H}_1^{orth}$ is $(q^m-1)(B_{m-1}+1)+1$.
Similar to case (i), the parameters of QLDPC code $Q$ with stabilizer $\textbf{S}$ follow immediately.

\textit{Remark.} (i) We need to gurantee that the parameters of QLDPC codes are meaningful, so we take the values of $q$ and $m$ into two cases.

(ii) We here consider the points and lines in $EG(m,q)$, where $q=p^s$. Only the case $s\geq 2$ were discussed in [6]. But here we can see that QLDPC codes also exist for case $s=1$, and we will give an example in the sequel.\\

\textbf{Example 2.} Let $m=q=2$. Consider the points without the origin and lines not passing through $\textbf{0}$ of $EG(2,2)$. The incidence matrix of points and lines is
$$
\textbf{H}_1=
\begin{pmatrix}
	0 & 1 & 1 \\
	1 & 0 & 1 \\
	1 & 1 & 0 \\
\end{pmatrix}.
$$
We can get
$$
\textbf{H}_1^{orth}=
\begin{pmatrix}
	\textbf{H}_1^T &\textbf{1} &\textbf{I} \\
\end{pmatrix},
$$
and
$$
\textbf{S}=
\begin{pmatrix}
	\textbf{H}_1^{orth} &\textbf{O} \\
    \textbf{O}  &\textbf{H}_1^{orth}\\
\end{pmatrix},
$$
The QLDPC code $Q$ with stabilizer $\textbf{S}$ has parameters $[[7,1,3]]$. $Q$ is actually the Steane code[15].\\

For the case $q$ odd and $m=2$, which is not mentioned in Theorem 3, it will make the dimension of $Q$ less than 0, if considered as the above. In the next part, we will give a new construction method for this case.

\subsection{A new construction of QLDPC code}
It shows that QLDPC codes cannot be constructed when $q$ is odd and $m=2$ in subsection A. In order to solve it, we give a new construction.  \\

\textbf{Theorem 4.}
Here $\textbf{H}_1$ is still the incidence matrix of points and lines of EG($m,q$) as in part A.
If $q$ is odd and $m=2$, then $\textbf{H}_1$ has size of $(q^2-1)\times(q^2-1)$. Let
$$
\textbf{H}_1^{orth}=
\begin{pmatrix}
	\textbf{H}_1^T &\textbf{1} &\textbf{I} &\textbf{I}\\
\end{pmatrix},
$$
and
$$
\textbf{S}=
	\begin{pmatrix}
		\textbf{H}_1^{orth} &\textbf{O}\\
		\textbf{O} &\textbf{H}_1^{orth}
	\end{pmatrix}.
$$
The QLDPC code $Q$ with stabilizer $\textbf{S}$ has parameters $[[3q^2-2, q^2, 2]]$.

\textit{Proof.} When $q$ is odd and $m=2$, the row weight of $\textbf{H}_1^T$ is $q$. Hence the row weight of $\textbf{H}_1^{orth}$ is $q+3$, which is an even number. Since any two distinct rows of $\textbf{H}_1^T$ has one ``1" in common, any two rows of $\textbf{H}_1^{orth}$ are orthogonal. Therefore, $\textbf{H}_1^{orth}$ is self-orthogonal.

The LDPC code having $\textbf{H}_1^{orth}$ as parity check matrix has minimum distance $2$, therefore, the minimum distance of $Q$ is also $2$.
The codelength and dimension of $Q$ follows immediately in the way similar to Theorem 3. \\

\textbf{Example 3.} Suppose $m=2, q=3$. The incidence matrix
$$
\textbf{H}_1=
\begin{pmatrix}
	0 & 1 & 0 & 0 & 1 & 0 & 1 & 0\\
	0 & 0 & 1 & 0 & 0 & 1 & 0 & 1\\
	1 & 0 & 0 & 0 & 1 & 0 & 0 & 1\\
    0 & 0 & 1 & 1 & 0 & 0 & 1 & 0\\
    1 & 0 & 0 & 0 & 0 & 1 & 1 & 0\\
    0 & 1 & 0 & 1 & 0 & 0 & 0 & 1\\
    1 & 1 & 1 & 0 & 0 & 0 & 0 & 0\\
    0 & 0 & 0 & 1 & 1 & 1 & 0 & 0\\
\end{pmatrix}.
$$
We have
$$
\textbf{H}_1^{orth}=
\begin{pmatrix}
	\textbf{H}_1^T &\textbf{1} &\textbf{I} &\textbf{I}\\
\end{pmatrix},
$$
and
$$
\textbf{S}=
	\begin{pmatrix}
		\textbf{H}_1^{orth} &\textbf{O}\\
		\textbf{O} &\textbf{H}_1^{orth}
	\end{pmatrix}.
$$
The QLDPC code $Q$ with stabilizer $\textbf{S}$ has parameters $[[25, 9, 2]]$.

\section{QLDPC CODES BASED ON ALL THE POINTS AND LINES OF $EG(m,q)$}
\subsection{All the points and lines}
In this part, we construct QLDPC codes using all the points and lines of $EG(m,q)$, where $m\geq2$.

Consider all the points and lines in $EG(m,q)$. Take the binary matrix $\textbf{H}_2=(h^{(2)}_{i,j})$ to be the incidence matrix of points and lines, whose rows are indexed by lines, and columns are indexed by points. $h^{(2)}_{i,j}=1$ if line $l_i$ passes through the point $v_j$, otherwise, $h^{(2)}_{i,j}=0$. Then $\textbf{H}_2$ has size of $(q^{m-1})(q^m-1)/(q-1) \times q^m$. There are $q$ ``1" in each row, and $(q^m-1)/(q-1)$ ``1" in each column. Since any two distinct points are joined by exactly one line, any two different columns of $\textbf{H}_2$ has just one ``1" in common. Two distinct rows of $\textbf{H}_2$ has at most one ``1" in common.

The row weight of $\textbf{H}_2^T$, $$\lambda=(q^{m}-1)/(q-1)=1+q+q^2+\cdots+q^{m-1},$$
is odd only when $q$ is even, or both $q$ and $m$ are odd. On the other hand, $\lambda$ is even if and only if that $q$ is odd and $m$ is even.
Now we define $\textbf{H}_2^{orth}$ as follows.
$$
\textbf{H}_2^{orth}=\left\{
	\begin{aligned}
		\;&
		\begin{split}
			[\;\textbf{H}_2^T\;|\;\textbf{1}\;],\;\ \ \ \ \ \ 
			\text{if $ q $ is even with $ m\geq3 $, or both $ q $ and $ m $ are odd.}
		\end{split} \\
		\;& [\;\textbf{H}_2^T\;|\;\textbf{1}\;|\;\textbf{I}\;],\;\ \ \ \text{if $ q $ is odd and $ m $ is even.}   \\
		\;& [\;\textbf{H}_2^T\;|\;\textbf{1}\;|\;\textbf{I}\;|\;\textbf{I}\;],\; \text{if $ q $ is even and $ m=2 $.}
	\end{aligned}
	\right.
	$$\\

\textbf{Lemma 5.} $\textbf{H}_2^{orth}$ is self-orthogonal.

\textit{Proof.} We need to show that $\textbf{H}_2^{orth} \cdot (\textbf{H}_2^{orth})^T=\textbf{O}$, that is, every row of $\textbf{H}_2^{orth}$ has even number of ones and any two distinct rows has even number of ones in common.

Since any two different columns of $\textbf{H}_2$ has just one ``1" in common, any two distinct rows of $\textbf{H}_2^{orth}$ intersect in exactly two nonzero positions.

(i) If $q$ is even and $m\geq3$, or both $q$ and $m$ are odd, the row weight of $\textbf{H}_2^{orth}$ is
$$B_m+1=2+q+\cdots+q^{m-1}.$$
Obviously, the row weight is an even number.

(ii) When $q$ is odd and $m$ is even, the row weight of $\textbf{H}_2^{orth}$
$$B_m+2=3+q+\cdots+q^{m-1}$$
is clearly even.

(iii) If $q$ is even and $m=2$, the row weight of $\textbf{H}_2^{orth}$ is
$$B_m+3=4+q+\cdots+q^{m-1},$$
which is an even number. \\

\textbf{Theorem 6.} Let
$$\textbf{S}=
	\begin{pmatrix}
		\textbf{H}_2^{orth} &\textbf{O}\\
		\textbf{O} &\textbf{H}_2^{orth}
	\end{pmatrix}.$$
The QLDPC code $Q$ with the stabilizer $\textbf{S}$ has parameters $[[n_2,k_2,d_2]]$ in the following cases.

(i) If $q$ is even and $ m\geq3 $, or both $ q $ and $ m $ are odd,
$$\begin{aligned}
& n_2=q^{m-1}B_m+1,\\
& k_2=q^{m-1}B_m-2q^m+1,\\
& d_2\geq q+1.
\end{aligned}
$$

(ii) If $ q $ is odd and $ m $ is even,
$$\begin{aligned}
& n_2=q^{m-1}B_{m}+q^m+1,\\
& k_2=q^{m-1}B_{m}-q^m+1,\\
& d_2\geq q+1.
\end{aligned}
$$

(iii) If $q$ is even and $ m=2 $,
$$\begin{aligned}
& n_2=3q^{2}+q+1, \\
& k_2=q^{2}+q+1,\\
& d_2=2.
\end{aligned}
$$

\textit{Proof.} (i) If $q$ is even and $m\geq3$, or both $q$ and $m$ are odd, then
$$
\begin{aligned}
	\textbf{S}=
	&\bordermatrix{
		&q^{m-1}B_{m}+1&q^{m-1}B_{m}+1\cr
		&\textbf{H}_2^{orth} & \textbf{O} \cr
		&\textbf{O} &\textbf{H}_2^{orth} \cr
	}\begin{array}{l}
		q^m\\q^m  \end{array}.
\end{aligned}
$$
The QLDPC code $Q$ with stabilizer $\textbf{S}$ is a CSS code.
The codelength of $Q$ is obviously $q^{m-1}B_{m}+1$. Dimension of $Q$ is $q^{m-1}B_m-2q^m+1,$ which is positive only when $m\geq3$.
The LDPC code having $\textbf{H}_1^{orth}$ as its parity check matrix has minimum distance at least $q+1,$ so it is with $Q$.

(ii) When $q$ is odd and $m$ is even,
$$
\begin{aligned}
	\textbf{S}=
	&\bordermatrix{
		&q^{m-1}B_{m}+q^m+1&q^{m-1}B_{m}+q^m+1\cr
		&\textbf{H}_2^{orth} & \textbf{O} \cr
		&\textbf{O} &\textbf{H}_2^{orth} \cr
	}\begin{array}{l}
		q^m\\q^m  \end{array}.
\end{aligned}
$$
Parameters of $Q$ are easy to verify in the way similar to case (i).

(iii) If $q$ is even and $m=2$,
$$
\begin{aligned}
	\textbf{S}=
	&\bordermatrix{
		&3q^{2}+q+1&3q^{2}+q+1\cr
		&\textbf{H}_2^{orth} & \textbf{O} \cr
		&\textbf{O} &\textbf{H}_2^{orth} \cr
	}\begin{array}{l}
		q^2\\q^2  \end{array}.
\end{aligned}
$$
The LDPC code having $\textbf{H}_1^{orth}$ as its parity check matrix has minimum distance $2,$ so it is with the minimum distance of Q.
Similar to case (i), the codelength and dimension of $Q$ follow immediately.

\textit{Remark.} We need to guarantee that the dimension of the QLDPC code $Q$ is positive, so we separate the values of $q$ and $m$ into three cases as above. \\

\textbf{Example 4.} Let $m=q=2$. There are 4 poins and 6 lines in $EG(2,2).$ The transpose of the incidence matrix of all the points and lines is
$$
\textbf{H}_2^T =
\begin{pmatrix}
	1 & 1 & 1 & 0 & 0 & 0 \\
	1 & 0 & 0 & 0 & 1 & 1 \\
	0 & 1 & 0 & 1 & 0 & 1 \\
	0 & 0 & 1 & 1 & 1 & 0
\end{pmatrix}.
$$
The row weight of $\textbf{H}_2^T$ is $3$, we can get
$$
\textbf{H}_2^{orth} =
\begin{pmatrix}
	\textbf{H}_2^T &\textbf{1} &\textbf{I} &\textbf{I}
\end{pmatrix},
$$
and
$$\textbf{S}=
	\begin{pmatrix}
		\textbf{H}_2^{orth} &\textbf{O}\\
		\textbf{O} &\textbf{H}_2^{orth}
	\end{pmatrix}.$$
The QLDPC code $Q$ with stabilizer $\textbf{S}$ has parameters  $[[15,7,2]]$.

\subsection{Parallel Classes}
In this part, we consider the QLDPC codes based on the parallel classes of lines in EG($m,q$).

Take into account all the points and lines of $EG(m,q)$. There are $q^{m-1}(q^m-1)/(q-1)$ lines in $EG(m,q)$. They can be partitioned into $(q^m-1)/(q-1)$ parallel classes, and each class contains $q^{m-1}$ lines. Lines in each parallel class have no points in common. We denote the classes to be $\mathcal{H}_1, \cdots, \mathcal{H}_t$, where $t=(q^m-1)/(q-1)$. For $i=1, \cdots, t$, take $\textbf{A}_i$ to be the incidence matrix of the points in $EG(m,q)$ and lines in $\mathcal{H}_i$, whose rows are indexed by lines and columns are indexed by points. Let $\textbf{A}_i^{orth}$ be as follows.
$$
\textbf{A}_i^{orth}=\left\{
\begin{aligned}
	& \;\textbf{A}_i,\;\ \ \ \ \ \ \ \ \ \ \text{if $q$ is even and $q\geq4$}.\\
	& \;[\;\textbf{A}_i\;|\;\textbf{I}\;],\;\ \ \ \text{if $ q $ is odd}.\\
	& \;[\;\textbf{A}_i\;|\;\textbf{I}\;|\;\textbf{I}\;],\;\text{if $q=2$}.
\end{aligned}
\right.
$$\\

\textbf{Lemma 7.} $\textbf{A}_i^{orth}$ is self-orthogonal for $i=1, \cdots, t$.

\textit{Proof.}
Rows of $\textbf{A}_i$ correspond to the lines in parallel class $\mathcal{H}_i$. Since lines in each parallel class have no points in common, any two distinct rows of $\textbf{A}_i$ are orthogonal. So it is with $\textbf{A}_i^{orth}$.

(i) If $q$ is even and $m\geq4$, the row weight of $\textbf{A}_i^{orth}$ is $q$.
Obviously, it is an even number.

(ii) When $q$ is odd, the row weight of $\textbf{A}_i^{orth}$ is $q+1$, which is even.

(iii) If $q=2$, the row weight of $\textbf{A}_i^{orth}$, $4$, is even. \\

\textbf{Theorem 8.} The QLDPC code $Q_i$ with the stabilizer $\textbf{S}_i$ has parameters $[[n_3,k_3,d_3]]$, where $i=1, \cdots, t$,
$$
\textbf{S}_i=
\begin{pmatrix}
	\;\textbf{A}_i^{orth} &\textbf{O} \\
	\textbf{O} &\textbf{A}_i^{orth}
\end{pmatrix},
$$
and
$(n_3,k_3,d_3)= $
$$
\left\{
\begin{aligned}
	\;& (q^m, q^m-2q^{m-1}, 2),\; \text{if $ q $ is even with $q\geq4$.} \\
	\;& (q^m+q^{m-1}, q^m-q^{m-1}, 2),\; \text{if $q$ is odd.}\\
	\;& (2^{m+1}, 2^m, 2),\; \text{if $ q =2 $.}
\end{aligned}
\right.
$$

\textit{Proof.} (i) If $q$ is even and $q\geq4$,
$$
\textbf{S}_i=
\bordermatrix{
	&q^m &q^m\cr
	&\textbf{A}_i^{orth} &\textbf{O} \cr
	&\textbf{O} &\textbf{A}_i^{orth} \cr
}\begin{array}{l}
	q^{m-1}\\q^{m-1}  \end{array}.
$$
The codelength $n_3$ of $Q_i$ is the number of columns of $\textbf{A}_i^{orth}$, which is $q^m$.
The dimension $$k_3=q^m-2q^{m-1}.$$ The LDPC code having $\textbf{A}_1^{orth}$ as parity check matrix has minimum distance $2,$ so it is with the minimum distance of $Q_i$.

(ii) If $q$ is odd,
$$
\textbf{S}_i=
\bordermatrix{
	&q^m+q^{m-1} &q^m+q^{m-1}\cr
	&\textbf{A}_i^{orth} &\textbf{O} \cr
	&\textbf{O} &\textbf{A}_i^{orth} \cr
}\begin{array}{l}
	q^{m-1}\\q^{m-1}  \end{array}.
$$
Similar to cases (i), the parameters of $Q_i$ follow immediately.

(iii) When $q=2$,
$$
\textbf{S}_i=
\bordermatrix{
	&2^{m+1} &2^{m+1}\cr
	&\textbf{A}_i^{orth} &\textbf{O} \cr
	&\textbf{O} &\textbf{A}_i^{orth} \cr
}\begin{array}{l}
	2^{m-1}\\2^{m-1}  \end{array}.
$$
Similar to cases (i), the parameters of $Q_i$ are easy to verify.\\

\textbf{Example 5.} Let $m=q=2$. Incidence matrices of points and lines of parallel classes  in $EG(2,2)$ are
$$
\textbf{A}_1=
\begin{pmatrix}
	1 & 1 & 0 & 0\\
	0 & 0 & 1 & 1
\end{pmatrix},\ 
\textbf{A}_2=
\begin{pmatrix}
	1 & 0 & 1 & 0 \\
	0 & 1 & 0 & 1
\end{pmatrix},\ 
\textbf{A}_3=
\begin{pmatrix}
	1 & 0 & 0 & 1  \\
	0 & 1 & 1 & 0
\end{pmatrix}.
$$
Then
$$
\textbf{A}_i^{orth}=
\begin{pmatrix}
	\textbf{A}_i &\textbf{I} &\textbf{I}
\end{pmatrix},\ 
\textbf{S}_i=
\begin{pmatrix}
	\textbf{A}_i^{orth} &\textbf{O} \\
	\textbf{O} &\textbf{A}_i^{orth}
\end{pmatrix}.
$$
The QLDPC code $Q_i$ with stabilizer $\textbf{S}_i$ has parameters $[[8,4,2]]$.

\section{Conclusion}
In this paper, we gave the requirements of QLDPC codes constructed in [6], and presented a new construction method for cases not included in the above conditions. Similarly, we constructed new QLDPC codes using all the lines and lines of parallel classes of Euclidean geometry, respectively. The constructed codes have high rates and their minimum distances are bounded.

\section*{Acknowledgments}
We would like to thank the editor and the anonymous reviewers for all their efforts in reviewing the paper and providing detailed comments.

\vfill

\end{document}